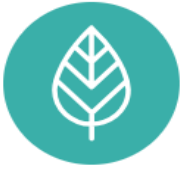

**Journal of Future Digital Optimization**

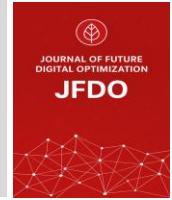

# Systematic Review of Smart Factories Production in Industry 5.0

## Ali Bakhshi Movahed[a], Hamed Nozari[b], Aminmasoud Bakhshi Movahed[c*]

*[a]School of Industrial Engineering, Iran University of Science and Technology, Tehran, Iran*

*[b]Department of Management, Azad University of the Emirates, Dubai, UAE*

*[c]School of Management, Economics and Progress Engineering, Iran University of Science and Technology, Tehran, Iran*

**Abstract**

Technology is undeniably positioned in today's industrial world, especially in manufacturing and smart factories. Unlike other industrial revolutions, humans are more important in the fifth generation of the Industrial Revolution. According to this fact, one of the critical factors of Industry 5.0 (I 5.0) can explained as human-centricity. Also, the implementation of modern technologies can be observed in smart factories. These types of factories produce more comfort and professionalism than the others. The current study states the importance of I 5.0 and innovative factory production (SFP). Furthermore, 36 articles are reviewed and, according to the meta-synthesis methodology, are categorized well. The research emphasizes the impact of I 5.0 on SFP by modern technologies and the whole policy. The new world can streamline people's new lives and create a revolution in the production lines of smart factories. Improving the structure of factories can be as feasible as this optimistic insight is.

*Keywords:* Industry 5.0; Smart factory; Meta-synthesis.

## 1. Introduction

Some critical policies in the industry 5.0 (I 5.0) generation are vital to know. According to the studies, Policies about I 5.0 mostly encourage the moral utilization of modern technologies and digital progress, facilitate the improvement of abilities and acquisition of new ones by the employees, and guarantee the fair allocation of advantages. Also, the stakeholder's insight into I 5.0 can utilize the technological and functional elements [1]. As a practical insight, assess the stakeholders' roles in the food supply chain in adopting I 5.0 drones to promote cleaner production [2]. This sample can be strong evidence to support that the food and beverage industry has some examples in I 5.0. Furthermore, policy papers frequently entail assistance from both the public and private sectors. Investment in R&D and the establishment of regulatory frameworks can be helpful. Also, I 5.0 has some impact on production. Industry and production have some close connections to each other. In addition, the fifth generation of the Industrial Revolution played a role in the whole production process. Factories are generally affected by the Industrial Revolution, and the fifth one can also develop and enhance smart factory production (SFP).

Some new technologies can be helpful in manufacturing, including automation, robotics, and AI, which can significantly improve operational topics in smart factories' efficiency. Using modern technologies, as stated before, can facilitate the development of smart factories. SFP is related to new technologies. All the production lines and all the assembly lines are affected by the factory's policy. Surely, if the policy uses more advanced levels of technology, the manufacturing systems and processes will be upgraded, too. It is happening that with new and modern production systems, the products will be new too. Also, by having new products, marketers of the companies and factories can easily sell the products. At the same time, the workers and employees will be more satisfied, and, with this atmosphere, the whole organizational culture will be upgraded, too. As a result, it can be observed that the whole organization can be affected by the production lines and systems. As it is in the definition of systems, all events in the company start with the manufacturing process. All these events happen just like a chain. So, it is necessary to state that the SFP in Industry 5.0 can precisely improve the supply chain of the factories.

By streamlining production, reducing bugs, and optimizing resource use, these technologies lead to faster cycles, better quality control, and increased overall efficiency. Prospectively, one of the impacts of I 5.0 is improving customization [3]. Also, I 5.0 enables suitable implementation of customization and personalization in the production process through advanced technologies. Paying attention to the customers individually is an essential factor in the personalization process. As always, SFP in the personalization process can be effective. In the smart factory, the customization process is better because it analyzes the data and makes valid information. It is essential that more precise measurements can create more accurate results. This argument can support the idea that the brighter the factory is, the more accurate the results can be. In the smart factory, analyzing the data tends to be meticulous, and it helps them find out customers' tastes. For example, in B2C companies, the company's favorite product can easily be found through data analysis. Then, the cycle starts, and the managers will know the situation. So, now it is the time to use the information and criteria gathered.

Also, by collaborating with machines in I 5.0, employees can increase productivity by assigning creative and critical thinking tasks to people and repetitive duties to machines, creating a collaborative environment that increases productivity. Managers must hire new people adept at handling cutting-edge technologies, data management and analysis, and machine programming [4]. Furthermore, I 5.0 requires a shift in workforce skills due to automation and AI. Additionally, safety is improved in the new industry versions according to these facts, and I 5.0 can be effective in the production process.

Smart, intelligent, or Industry 4.0 (I 4.0) factories are advanced manufacturing facilities that apply cutting-edge technology and automation to improve production and efficiency. Smart factories efficiently accommodate market needs and demands, promote collaboration, and create dynamic virtual organizations [5]. These factories integrate IoT devices, AI, and cloud computing to create a connected and intelligent ecosystem. Also, Smart factories are industrial facilities that leverage digital technologies to enhance productivity [6]. I 4.0 factories can encompass advanced manufacturing facilities. These manufacturing facilities integrate digital and modern technologies and automation to optimize production systems and processes. These factories use a wide and innovative range of technologies, such as robotics and cloud computing, to establish intelligent and interconnected systems. Within I 4.0 factories, machines are equipped with sensors and connected to a network, enabling the collection and analysis of real-time data. The use of sensors in the SFP has an attractive story, too. The data is then utilized to improve and optimize the decision-making process. This process starts with the mind and neurons. Moreover, smart factories mostly employ autonomous robots and intelligent systems to carry out the duties and tasks traditionally performed by humans. The first goal of I 4.0 factories is to establish highly adaptable, efficient, and responsive manufacturing systems that respond to changing demands and optimize resource use. By integrating modern and digital technologies into the production process, these factories strive to achieve levels of automation and customization. The I 4.0 industry policies also emphasize personalization and automation. This happens because these titles are the needs of people in our new world.

There are some examples of using technologies in the manufacturing process. IoT integration in smart factories can connect devices, sensors, and machines through networks. According to the speed of the internet, it is possible to apply. The utilization of advanced industrial Internet of Things (IIoT) and cyber-physical system (CPS) technologies has led to an essential shift in the method of assigning manufacturing resources [7]. Smart Factory Production can quickly help produce outputs.

Finding proper balance is the key to success. On one side, Machines can manage repetitive tasks efficiently, while on the other side, humans generate innovation and creativity. Also, problem-solving skills and adaptability to complex conditions are advantages of that. Furthermore, the advantages and disadvantages are concluded in this paper. Collaboration between the two often yields the best results in a business setting. Companies need to be controlled by machines or humans. Research shows that most repetitive factory processes do not need human monitoring; machines do, however, to manage the situation. For this reason, monitoring by machine is a favorite topic for many factories. Also, remote control of machines can be an attractive subject for company managers. Furthermore, this subject can facilitate the burden on the worker's shoulders in the factories. Remote monitoring of machines is an important factor in enhancing visibility and managing the situations of the SFP [8]. Also, predictive maintenance and intelligent sensors in smart factories can improve the whole maintenance process in SFP [9]. This topic has been seriously critical in recent years. Sensors, with the collaboration of AI engineers, enable the top management of factories to make better decisions. Making balance in using inputs and outputs can easily be performed soon.

The concept of Time can update and outdate many things in the world. Technology and time are the two concepts that are connected. New technologies can update intelligent systems and production lines. As a practical insight, some smart sensors can shift the old paradigm to the new one in the Production lines of different factories [10]. Furthermore, production planning and strategic scheduling exist for all factories, especially the smart ones [11]. Smart factories use various products to automate and streamline manufacturing. Also, advanced technologies can produce better outputs than ever in smart factories.

I 5.0, the human-centric industrial revolution, combines human abilities with innovative technologies in smart factories. This approach emphasizes collaboration between humans and machines to earn better productivity and innovation. It has a significant impact on SFP. Thus, the factories collaborate between humans and technologies [12]. One of the most critical effects of I 5.0 on SFP is increased communication efficiency between people from an industrial background. This means that I 5.0 advocates for improved collaboration between humans and machines. This mutual collaborative approach has odds of making the optimum decision and increasing compatibility. The optimum decision is the one that can consume resources precisely and earn better results. It determines the balance between outputs and inputs. Also, the relationship between humans and machines has a long story that can be investigated in the future [13]. Another critical effect of I 5.0 on SFP is worker empowerment. It happens through the provision of essential tools, thereby enhancing their efficiency an d job satisfaction. These technologies may encompass wearable devices, augmented reality (AR) tools, and cobots that help workers and employees execute their complex tasks and ultimately enhance the overall performance of the workers in their environment. Cobots will be explained in the literature review in the following.

For example, there is an application of AR in a lean workplace that observes case studies at smart factories [14]. In addition, advancing modern technology can be another impact of I 5.0 on SFP. By embracing the advancements, smart factories can seamlessly cater to the unique needs of individual customers, efficiently manufacturing small batches or customized items. Also, a Smart system in Manufacturing with Mass Personalization (S-MMP) for Production lines and scenarios of production systems can be helpful in this field [15]. Thus, it shows that personalization and customization are enhanced. It has not only improved customer satisfaction but also improved market competitiveness. Also, progressing workforces can be another effect of I 5.0 on SFP. Promoting continuous learning and development among workers and employees to proficiently operate, maintain, and optimize the advanced technologies applied in innovative factory production is crucial. Also, I 5.0 makes it possible for factories. Additionally, I 5.0 strongly emphasizes ethical and sustainable manufacturing practices. AI and IoT can grow SFP and effectively monitor and optimize energy consumption.

Indeed, the issue of energy and its provision poses a significant challenge in today's world. Nowadays, the issue of energy and providing it is a challenging topic. Due to the increasing global population and the surging demand for energy, the provision of that will be more controversial. The feasibility, cost-effectiveness, and environmental consequences of various energy sources are being explored and debated at international conferences more and more. In the future, one of the most challenging concerns will be supplying energy for many countries. The gathered articles state the significance of energy. In general, I 5.0 can transform SFP by prioritizing humans in technological improvement. The aim is to establish a manufacturing environment that is more cooperative, productive, and environmentally friendly. It will benefit both employees and enterprises to collaborate more. This paper defines the structure as shown in the introduction in the first part. The introduction part defines and highlights the advantages of SFP and the importance of this topic. Then, it is time to provide the literature review. Also, some papers investigate the following. The methodology section consists of the Meta-synthesis and its charts. Furthermore, the tables are defined in the finding part, and the enablers and drivers are explained precisely. Finally, the future directions and the current trend are defined at the conclusion. The research steps present the sequences of events level by level to form the current paper. At the first layer, the review of SFP is a valuable step because the reader can understand how SFP is in the literature. Also, in the second stage, the meta-synthesis methodology is shown for the drivers and barriers in the background of SFP. Understanding the obstacles and how they are prioritized in the paper is essential. This happens for enablers, too, because they are crucial for readers to know. In the third layer, the analysis of barriers and enablers is presented. It is essential to analyze them because by analyzing them, people can easily understand their effects. The next level of the research is to present a model. A model can show any information that is the key to understanding the points. A conceptual model is presented in the background of SFP in I 5.0. In the last layer, the descriptive statistics of research and different types of taxonomies are shown. Also, the graphical abstract of the article is in Figure 1.

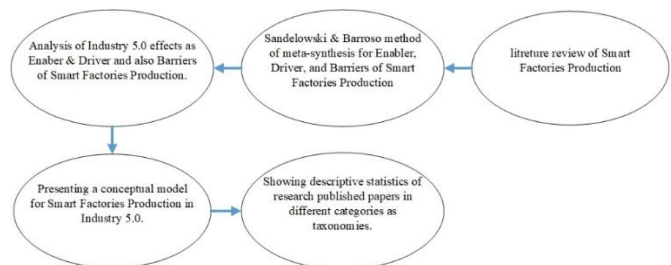

**Figure 1. Graphical Abstract.**

## 2. Literature review

Smart factories encompass the incorporation of cutting-edge technologies and automation to enhance the efficiency of manufacturing operations and processes. Smart factories heavily depend on interconnected devices, machines, and systems in modern factories. In the I 4.0 factory, machines and equipment are interconnected, communicate, and have a central control system. This type of connectivity enables real-time monitoring, analysis, and decision-making [16]. This connectivity facilitates smooth data exchange among components, enabling real-time monitoring and control.

The IIoT is instrumental in functioning intelligent factories as it facilitates the connection of physical devices and sensors to gather and transmit data. Integrating intelligent computing and AI paradigms in smart factories enables IIoT to streamline production processes and minimize the need for human interventions [17]. Also, this data can be significant as it offers a valuable point of view into machine performance. The data management process is trying to offer this.

According to the papers, some researchers claim that they can present a base for building structure and some critical success indicators to evaluate the capabilities of SFP [18]. Introducing these elements is vital for factories and companies because they can estimate their performance rate. The machine will be evaluated more precisely by the right indicators. Also, some smart factories use automation and robotics to optimize the production processes. Furthermore, Smart factories use big data analytics (BDA) and AI to predict and prevent equipment failures. So, smart factories must be able to analyze data and predict small and big failures to succeed on their road [19]. Reviewing articles shows the authors that being aware of successful people and paying attention to their lifestyles can help people understand the critical moments of their lives. Their routine is a complete course of lessons for people who want to be successful like them. The decisions they must make in challenging situations show the difficulties they face on their way to success.

Additionally, smart factories continuously monitor machine performance, detect irregularities, and schedule maintenance activities. Thus, this activity makes it easier for companies and factories to produce their products. As a result, the maintenance cost will be decreased after implementing that. Some technologies, like Digital Twins, are new in emerging countries' economies [20]. The role of Digital twins is to duplicate products or services virtually and use them in intelligent manufacturing facilities to enhance the production systems. It enables companies to test their modern ideas and identify the potential obstacles before implementing them. Also, Smart factories use systems that monitor in real-time and control the quality parameters. It also checks resource consumption and detects problems quickly.

The Industrial Revolution started with Industry 1.0 and continued with Industry 5.0. The fifth generation of the industry includes three elements. Human-centricity, resiliency, and sustainability are these three elements. However, I 5.0 impacts SFP by using worker robots with human workers to improve the efficiency and safety of whole organizations. The collaboration between these two types of workers will be excellent. If the factories want long-lasting collaboration, they should be aware of the new ways of maintenance.

Cobots are known as collaborative robots, too. Cobots are specific robotic systems that have progressed in operating with humans and workers. In contrast to conventional industrial robots, which are generally confined within cages or limited spaces, cobots have been specifically designed to ensure safety and directly interact with humans without posing a significant threat of harm. In a short explanation, a cobot is a softer version of a robot and can act more optimally. These robots are equipped with advanced sensors and programming capabilities that enable them to detect and respond to the presence of humans. Cobots can find extensive applications in various industries, such as manufacturing and assembly.

Collaborative robots called Cobots can help the production line of smart factories. Cobots can control the quality of the whole production system. Also, it can help humans handle logistics. Supply lines and logistics from the past until now are significant in all fields. It is good to know that the history of these concepts began with the wars. Cobots can monitor people's safety, too. Sensors and cameras have enabled the cobots to control the safety of workers in the factories. Cobots and I 5.0 can collaborate and enhance how people work in the firms and companies to begin a new generation in the factories.

This research collects 36 papers and categorizes them into several categories. Table 1 shows relevant studies from 2019 to 2023 in the last five years. Also, it can be evidence of innovative factory production.

**Table 1. Relevant papers with SFP.**

| Authors | Source Title | Methodology | Tech Type | Field | Application | Publisher |
|---|---|---|---|---|---|---|
| Zakoldaev et al. [21] | The multi-agent environment of cyber and physical for the I 4.0 *SFP* | Modelling | Cyber-Physical System | Other | Different Companies | IOP |
| Otles & Sakalli [22] | *I 4.0*: The *smart factory* of the future in the beverage industry | Mixed methodology | AI & IoT | Manufacturing | F&B industry | Elsevier |
| Zakoldaev et al. [23] | The life cycle of technical documentation in the *smart factory* of *I 4.0* | Modelling | I 4.0 | Digitalization | Different Companies | IOP |
| Osterrieder et al. [24] | The *smart factory* as a critical construct of *I 4.0* | Systematic Review | I 4.0 | Futurology | Intelligent System | Elsevier |
| Longo et al. [25] | *Industry 5.0:* A human-centric Perspective for Future Factory | Quantitative methodology | I 5.0 | Digitalization | Different Companies | MDPI |
| Yu et al. [26] | *SFP* & Operation Management | Quantitative methodology | I 4.0 | Management | Counseling | IEEE |
| Kalsoom et al. [27] | Advances in sensor technologies in the era of *smart factories* and *I 4.0* | Qualitative methodology | AI & IoT | Management | Counseling | MDPI |
| Büchi et al. [28] | *Brilliant factory* | Qualitative methodology | I 4.0 | Management | Counseling | Elsevier |

| Authors | Source Title | Methodology | Tech Type | Field | Application | Publisher |
|---|---|---|---|---|---|---|
| | performance and *I 4.0* | | | | | |
| Shi et al. [29] | *Smart factory in I 4.0* | Other | I 4.0 | Other | Intelligent System | Other |
| Setiawan et al. [30] | Multi-agent System in *Smart Factory* of *I 4.0* | Quantitative methodology | I 4.0 | Manufacturing | Intelligent System | Other |
| Képešiová & Kozák [31] | *Smart Factory in I 4.0* | Other | Other | Manufacturing | Intelligent System | Other |
| Mohammad et al. [32] | *Smart factory reference model for training on I 4.0* | Quantitative methodology | AI & IoT | Manufacturing | Production System | Other |
| Hawkins [33] | *Smart factory performance in I 4.0*-based manufacturing systems | Qualitative methodology | AI & IoT | Estimation | Production System | Other |
| Boniotti et al. [34] | A conceptual model for *SFP* data | Modelling | Big data | Management | Different Companies | Springer |
| Diering & Kacprzak [35] | Selected product on the SFP line | Quantitative methodology | Other | Manufacturing | Production system | Springer |
| Abyshev & Yablochnikov [36] | R&D of a Service-oriented Architecture for a *SFP* System | Quantitative methodology | Simulation | Architecture | Production System | Other |
| Margherita & Braccini [37] | The *socially sustainable factory* of Operator 4.0 and *Industry 5.0* | Systematic Review | I 5.0 | Futurology | Intelligent System | Other |
| Wang et al. [38] | A Digital Twin Model of *SFP* System | Other | Simulation | Digitalization | Energy | IOP |
| Rahman et al. [39] | Modeling on the *SFP* optimization | Quantitative methodology | Other | Estimation | Energy | Other |
| Schou et al. [40] | Deconstructing *I 4.0* with defining the *smart factory* | Qualitative methodology | I 4.0 | Digitalization | Different Companies | Springer |
| Chumnumporn et al. [41] | Effect of a Customer's *Smart Factory* Investment on a | Quantitative methodology | I 4.0 | Manufacturing | Different Companies | Other |

| Authors | Source Title | Methodology | Tech Type | Field | Application | Publisher |
|---|---|---|---|---|---|---|
| | Firm's *I 4.0* Technology Adoption | | | | | |
| THURAPAN & Pooripakdee [42] | *Smart Factory* with Lean Automation for a High-Performance and Sustainable Organization | Mixed methodology | Other | Management | F&B industry | University |
| Gupta & Randhawa [43] | Implementing *I 4.0* and Sustainable Manufacturing Leading to *Smart Factory* | Quantitative methodology | I 4.0 | Manufacturing | Energy | Springer |
| Tsao et al. [44] | *Smart factory* effect on companies in *I 4.0* | Mixed methodology | I 4.0 | Management | Energy | Other |
| Sreedhanya & Balan [45] | *I 4.0* Framework Using 7-layer Architecture for *Smart Factory* Application | Modelling | Big data | Architecture | F&B industry | Springer |
| Kajbaje et al. [46] | The New Age of Manufacturing-*I 4.0* and *Smart Factory* | Systematic Review | AI & IoT | Manufacturing | Different Companies | University |
| Lee et al. [47] | Key enabling technologies for *smart factories* in the automotive industry | Systematic Review | AI & IoT | Manufacturing | A Case study | Springer |
| Khang et al. [48] | Integration of Robotics and IoT to *Smart Factory* Infrastructure in *I 4.0* | Other | AI & IoT | Futurology | Other | Other |
| Magnus [49] | *Intelligent factory* mapping and design | Mixed methodology | I 4.0 | Futurology | Different Companies | Springer |
| Kumar et al. [50] | Review on prognostics and health | Systematic Review | AI & IoT | Management | Other | Elsevier |

| Authors | Source Title | Methodology | Tech Type | Field | Application | Publisher |
|---|---|---|---|---|---|---|
|  | management in smart factory |  |  |  |  |  |
| Ryalat et al. [51] | *Clever factory design based on cyber-physical systems and IoT towards I 4.0* | Mixed Methodology | Cyber-Physical System | Digitalization | A Case study | MDPI |
| Fortoul-Diaz et al. [52] | *Innovative Factory Architecture Based on I 4.0 Technologies* | Mixed Methodology | Cyber-Physical System | Architecture | A Case study | IEEE |
| Hartmann et al. [53] | *Industry 5.0 skills in a learning factory using existing technologies* | Mixed Methodology | I 5.0 | Education | Different Companies | SSRN |
| Potthoff & Gunnemann [54] | *Lean Learning Factory Integrating Aspects of Industry 5.0 and Sustainability* | Mixed Methodology | I 5.0 | Education | A Case study | SSRN |
| Rantschl et al. [55] | Extension of the *LEAD Factory* to address *Industry 5.0* | Quantitative Methodology | I 5.0 | Manufacturing | A Case study | SSRN |
| Cotta et al. [56] | Cognitive *Factory* in *Industry 5.0* | Mixed Methodology | I 5.0 | Other | Other | MDPI |
| Current Research | *SFP* in *Industry 5.0* | Systematic Review | I 5.0 | - | - | - |

The study contains questions that have been stated as follows:

- What are the drivers and enablers of SFP in I 5.0?
- What are the barriers to SFP in I 5.0?

## 3. Methodology

The data was analyzed in three sections using the meta-synthesis methodology. Also, there is a systematic study examining past studies. With meta-synthesis, novel insights are not taken in the original texts [57]. Summarized information has been obtained from other papers with the same subject [58] in meta-synthesis. This is an umbrella for developing new knowledge [59] based on the detailed and comprehensive analysis of qualitative research findings. The current study applies the seven-step method (Figure 2). In Figure 3, the review process is signified based on the inclusion and exclusion criteria. The CASP (Critical Skills Appraisal Program) criterion was applied to assess validity. The Kappa coefficient was used to evaluate the reliability, and its value was computed as 0.692, which is exactly higher than the acceptable range.

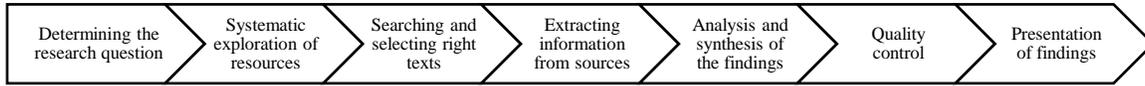

**Figure 2. Sandelowski & Barroso's method of meta-synthesis.**

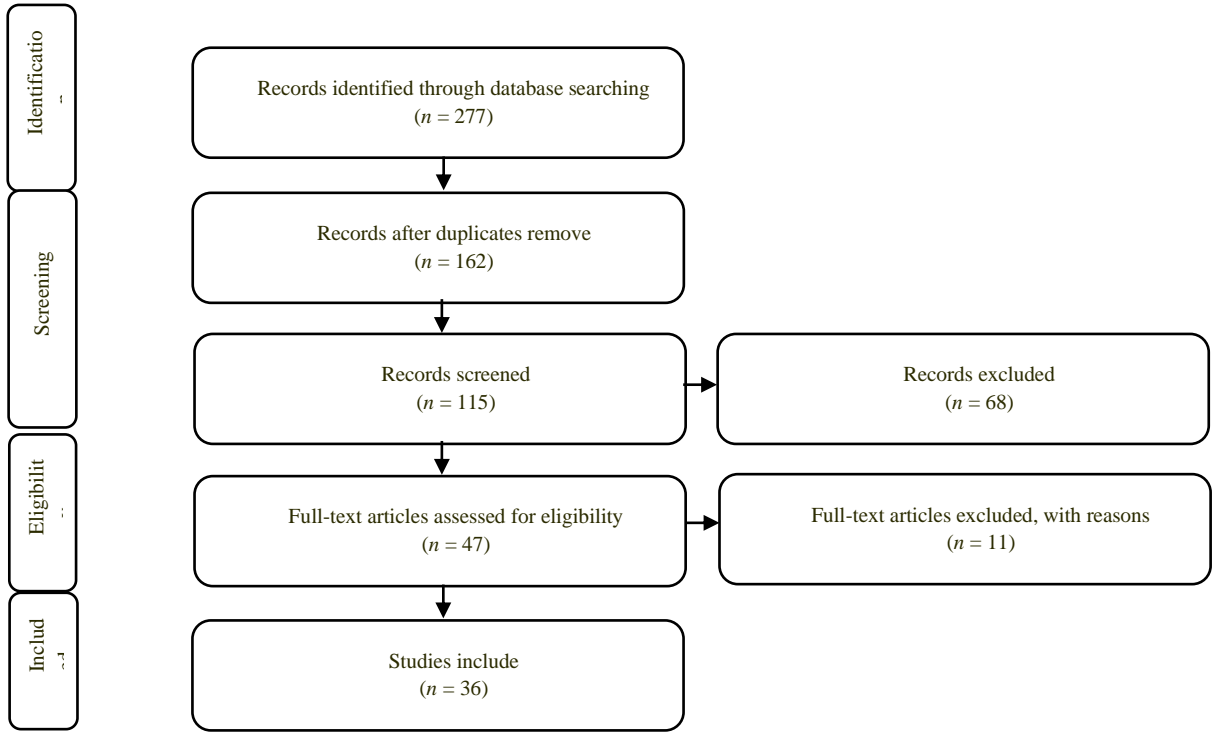

**Figure 3. PRISMA flow diagram of literature search.**

## 4. Finding

The meta-synthesis process has been done in two parts. The first section refers to the enablers. Drivers or enablers are arguments that make factories lean toward innovative productions in I 5.0. Enablers can firmly support improving the SFP. Drivers can precisely overcome obstacles in different factories [60]. Enablers with immediate changes [61] in value generation can change I 5.0 policies into novel circumstances [62]. The meta-synthesis analysis enablers are in Table 2.

Factories are attempting to leverage opportunities [63] for producing innovative products in new markets due to many barriers in the era of I 5.0 [64]. Barriers are resistant factors to changing productions toward SFP. Table 3 shows the analysis of barriers. The SFP framework in I 5.0 is shown in Figure 4.

### Table 2. Meta-synthesis of drivers and enablers

| Reference | First-order themes | Second-order themes | Aggregate dimensions |
|---|---|---|---|
| [44], [22] | Human-robot collaboration | Human detection prediction | I 5.0 Effects |
| [24], [38] | Human-Machine Cooperation | | |
| [24], [21], [60] | Human-Machine Interaction | | |
| [56], [45], [60] | Human-Machine Interface | | |
| [64], [39] | Evolving modern management ways | Sustainable innovation in automation | |
| [60], [56] | Production and process automation systems | | |
| [60], [38], [39] | Increased computing and production capacity | Adaptive capacity development with conscious resilience | |
| [43], [47] | A reliable source of resilience | | |
| [38], [47] | Cognitive bots for implementing repetitive processes with minimum cost and high accuracy | Optimizing Facilities | Production chain Transformation |
| [39], [44] | Advanced automation for | | |
| [47], [60], [52] | Digitization beyond automation with predictive analytics | Increasing responsiveness and predictability of the production line | |
| [60], [34] | Remote supervision and monitoring predictively | | |
| [53], [44] | Increasing productivity through entrepreneurship and time tracking | | |
| [22], [43] | Additive manufacturing for rapid production of prototypes with a limited volume of spare parts | Automation speed in setting up systems | |
| [38], [62] | Adopting automation system and co-automation | | |
| [21], [60] | Cyber-Physical Systems (CPS) can enhance safety | Safety, maintenance, and reliability considerations | |
| [60], [34] | maintenance is helped by Virtual Reality and Augmented Reality | | |
| [48], [52] | Reliability, validity, and scalability should be considered | | |
| [62], [35] | Data transmission from sensor to cloud | Sensor-based replenishment | Smart Supply |
| [60], [36], [21] | The connection between sensors and devices | | |
| [27], [21] | Updating sensors and Cyber-Physical Systems (CPS) | | |
| [48], [44] | Efficient delivery of data | Support and dynamic delivery | |
| [39], [27] | Transparent delivery process | | |
| [38], [21] | Autonomous robots to carry out warehouse operations | Digital supply network reinforcement | |
| [44], [34] | Using Robots Operation System (ROS) | | |
| [53], [60] | Cobots and Robots collaborate to improve safety systems | | |
| [48], [21] | A deep learning algorithm for the IIoT | Industrial Internet of Things (IIoT) | I 4.0 Technologies |
| [44], [60] | IIoT platforms | | |
| [35], [52], [60], [39] | Robots and AI integration | Artificial Intelligence (AI) | |
| [48], [60] | Adaption of AI-based systems and 3D printers | | |
| [53], [38] | Cobots and AI Collaboration | | |
| [65], [22] | Control agent activity based on AI and closed-loop control | | |
| [48], [53] | Big data analysis collaboration with AI | Big Data Analysis | |
| [21], [48] | Big data analysis collaboration with Machine learning algorithms and prediction | | |
| [60], [64] | Big Data and human-machine interaction and deep learning | | |
| [48], [60] | Cloud computing layer for connecting external systems | Cloud Computing | |
| [53], [43] | Cloud manufacturing | | |
| [21], [60] | Integrated energy data management | Control of operational errors and fluctuations | Data management by processes algorithmization |
| [22], [48] | Big data management systems | | |
| [48], [27] | increasing sustainability and efficiency of the system's operations by prediction | Operational inefficiencies prediction | |
| [24], [53] | manufacturing operations management (MOM) and operation planning | | |
| [39], [21] | optimized dynamic production systems | Self-optimize | Moving towards automating new governance models |
| [38], [64] | optimizes the supply chain due to warehouses | | |
| [48], [44] | Quick adaption to change market demands | Self-adapt | |
| [22], [46] | Self-adaptive decision-making process | | |
| [52], [46] | Role of a leader through a transformational leadership style | New career roles definition | Systems Specialists Skill Improvement |
| [64], [34] | Role of Industrial Systems Designers and human in managing positions | | |
| [45], [34] | Social skills in industrial environments and avoiding toxic workplace culture | Industrial culture development | Supportive beliefs of environmental society |
| [24], [34] | Culture and network infrastructure and a healthy work environment | | |
| [48], [47] | Adaption while changing skill obligations in the labor market | Professional obligations | |
| [22], [60] | Satisfying the requirements of the automotive industry | | |

### Table 3. Meta-synthesis of barriers.

| Reference | First-order themes | Second-order themes | Aggregate dimensions |
|---|---|---|---|
| [46], [47] | Risk for assembly and manufacturing errors | Error in human analysis and prediction | Forecasts and decisions monitoring |
| [22], [23] | Errors in prediction and detection in the process of the robot tasks | | |

| Reference | First-order themes | Second-order themes | Aggregate dimensions |
|---|---|---|---|
| [21], [52] | Synchronization of fundamental tools with the digital twins | Synchronization of senior decisions | |
| [24], [53] | Unable to provide a prediction of the consequences of decisions | | |
| [52], [44] | Incompatibility in the implementation of hardware | Extreme Hardware requirements | Challenges of optimizing operations |
| [60], [62] | Set up and holding costs of hardware and infrastructure | | |
| [24], [44] | Unavailable the license software | Extreme Software requirements | |
| [53], [44] | System Maintenance and its costs for software | | |
| [21], [36] | Unable to Provide a realistic cyber-physical environment | Environmental concerns | Lack of environmental resilience |
| [35], [34] | Costly physical damage and environmental effects | | |
| [35], [53] | Resistance to automation within the factory | Resistance to changes | |
| [22], [41] | Resistance to Production process change management | | |
| [22], [39] | Costs of Intelligent Manufacturing Integration System | Equipment upgrading cost | Huge investment requiring |
| [39], [35] | Generating data from automated equipment in a virtual environment insufficiently | | |
| [52], [41] | Managing employees with new technologies insufficiently | Employees retraining cost | |
| [45], [34] | Role of Older Employees and lack of cost-benefit analyses | | |
| [22], [38] | Lack of adoption in modular systems | Safe systems setting up cost | |
| [56], [28] | Lack of adoption in Prototype set up | | |
| [48], [47], [21] | Lack use of Cyber-physical human-centered system (CPHS) | Insufficient Cyber security | Security & Privacy |
| [24], [65] | Insufficient Industrial Cyber-Physical Platforms | | |
| [39], [48] | Privacy protection errors | privacy controls | |
| [43], [34] | Vulnerabilities to cyberattacks and privacy concerns | | |
| [47], [21] | insufficient use of dynamic capabilities and unskillfulness in management | Power conflicts in management | Structural situation instability |
| [39], [48] | Robots' engagement in workflows of management | | |
| [47], [21] | Lack of use of supply chain monitoring | Supply chain pressures | |
| [24], [38] | Weakness in the exchange of data by supply chain and product lifecycle | | |
| [52], [28] | Realizing the difficulties in organizational missions and targets | Lack of commitment | Cultural Issues |
| [38], [52] | Increasing costs and commitments by innovation studies | | |
| [52], [34] | Trustworthiness and trust are extremely | Lack of trust | |
| [34], [48] | The need for digital trust is extremely | | |
| [38], [43] | Lack of enough infrastructure in emerging technology | Inability to adapt to smart market changes | Alignment of the amortization process with the business model |
| [34], [64] | Incapable of responding to smart market fluctuations | | |
| [64], [48] | Lack of enough Integration of Product design and services | Transition to service business models; Resources & Machines as a service | |
| [46], [56] | Defective Physical product or service | | |

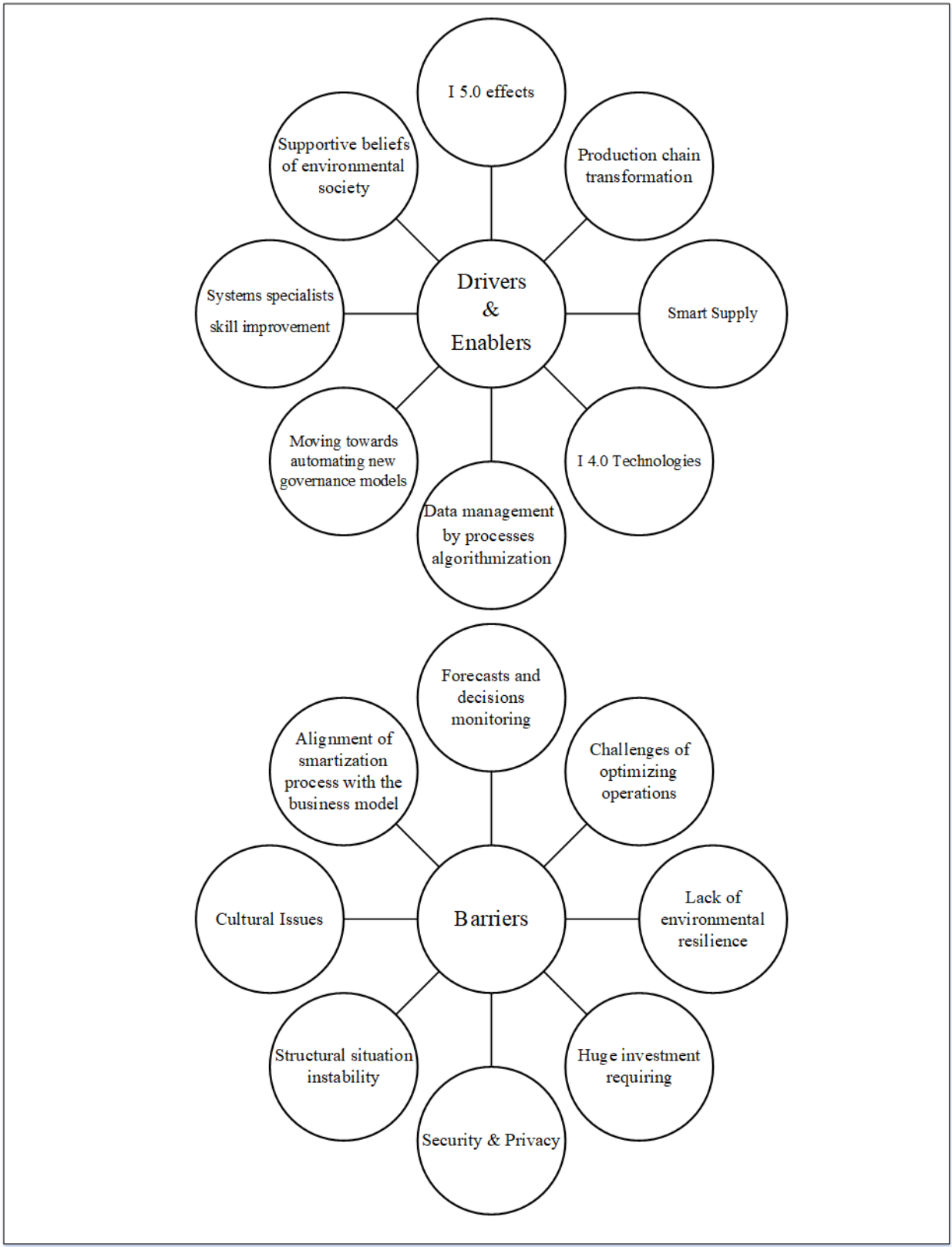

**Figure 4. SFP framework in I 5.0.**

### 4.1. Descriptive statistics of published papers

To understand the theme of research references, which has been the basis of the literature, different categories are presented in the following figures in this section. According to a search of publishers' databases from 2019 to 2023, 36 articles were found. Figure 5 shows the number of publications in different years. Figure 6 shows the number of publications by different publishers. Figure 7 shows an analysis of the field of the reviewed articles based on the percentage.

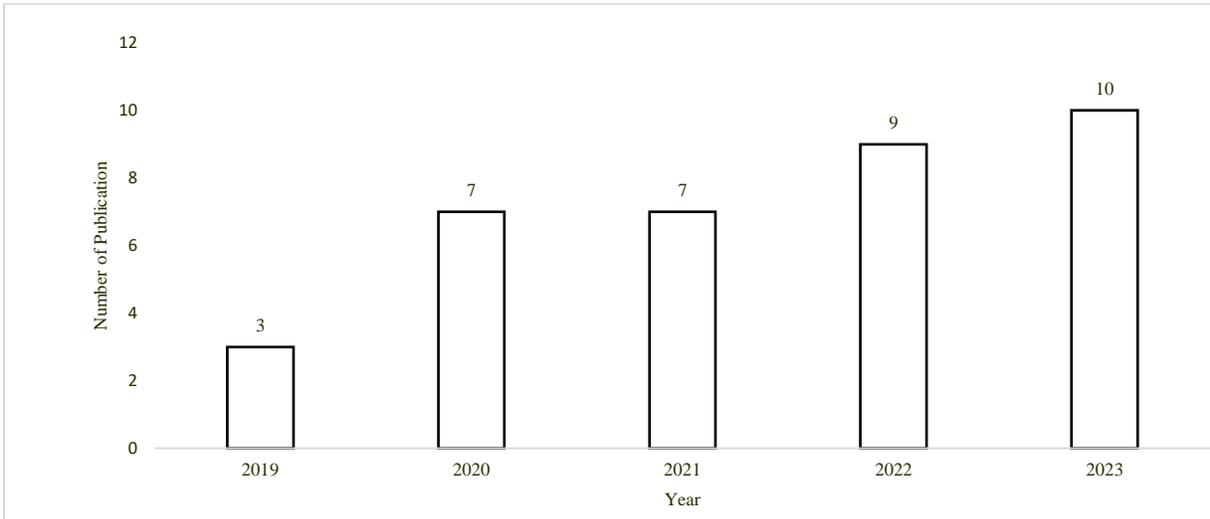

**Figure 5. Number of published articles in the last five years.**

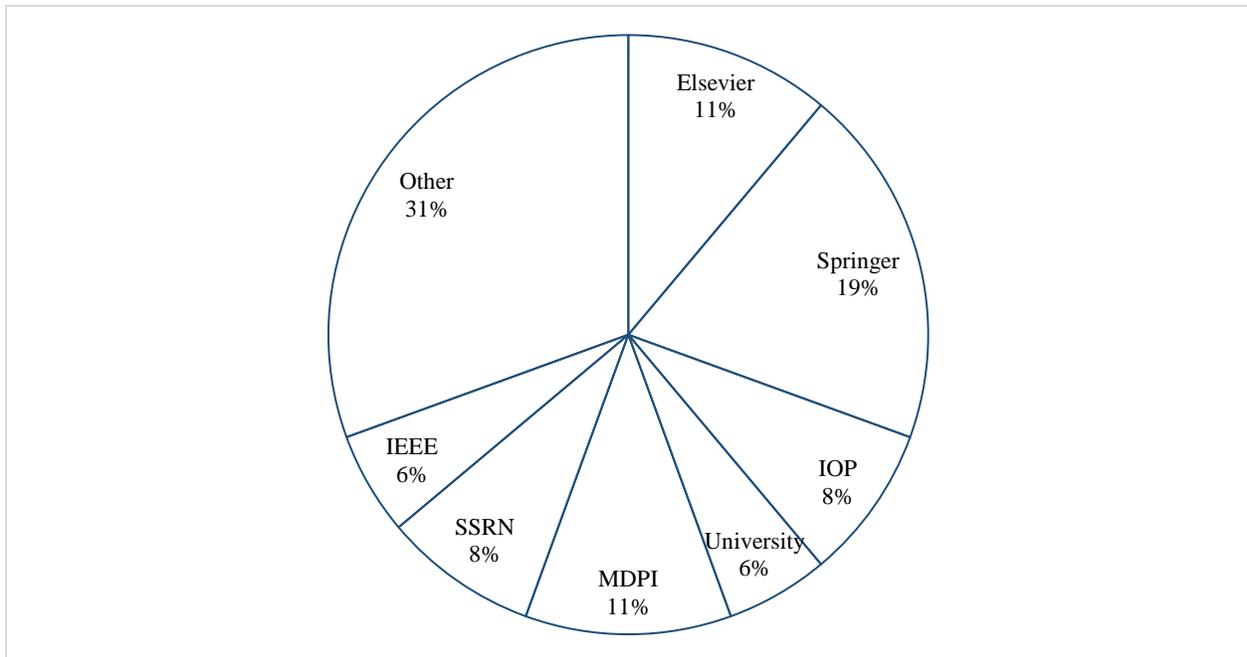

**Figure 6. Number of articles published by different publishers.**

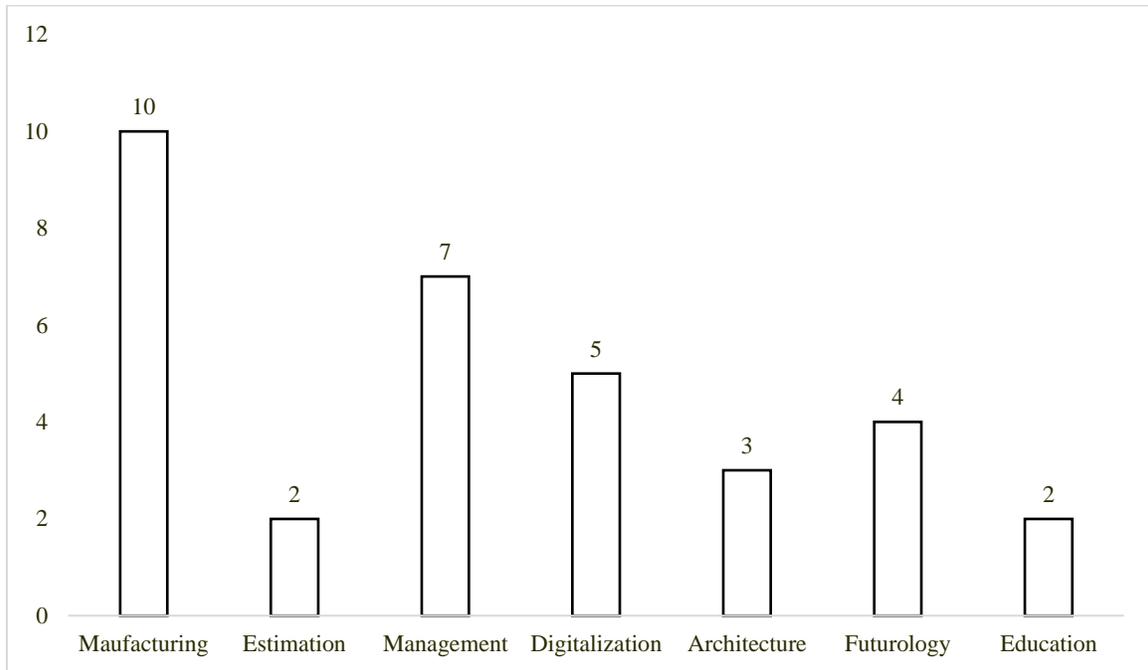

**Figure 7. Field of the articles.**

Figure 8 shows the application of the reviewed articles separately. Figure 9 shows the distribution of methods and data analysis of the articles. Figure 10 shows the segmentation using I 4.0 technologies.

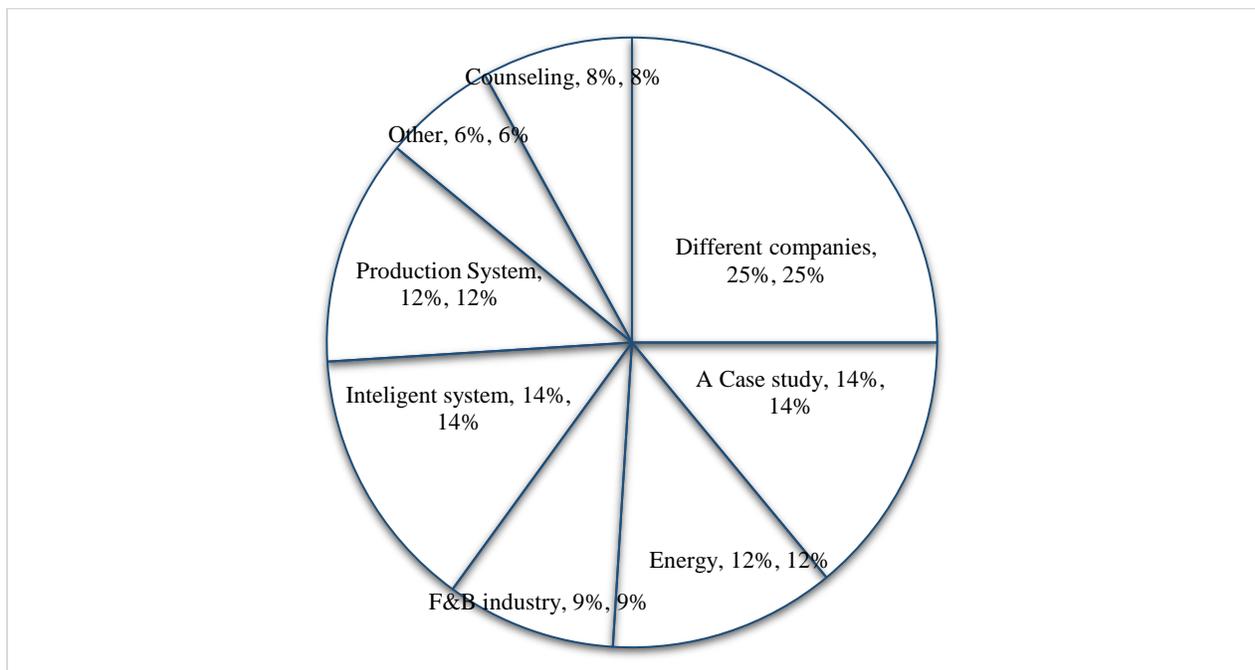

**Figure 8. Application of the articles.**

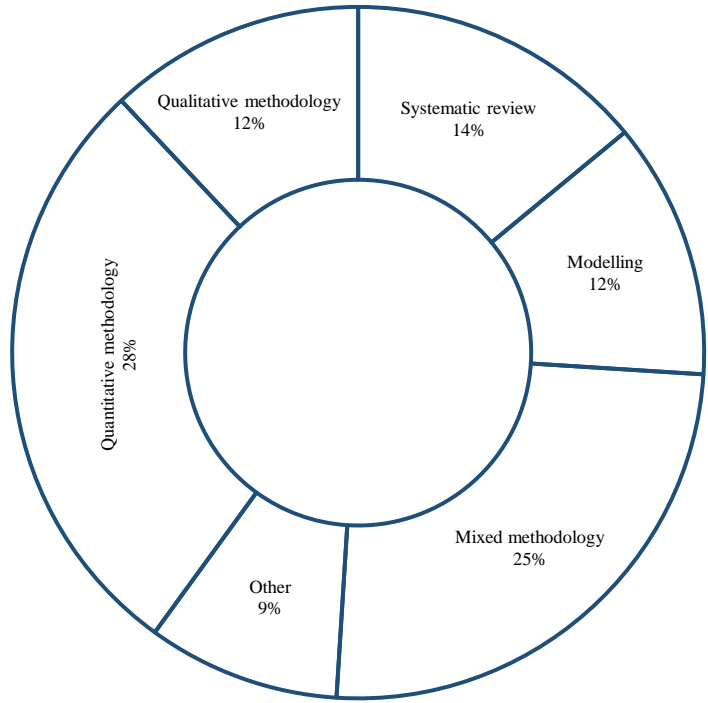

**Figure 9. Methods and data analysis of the articles.**

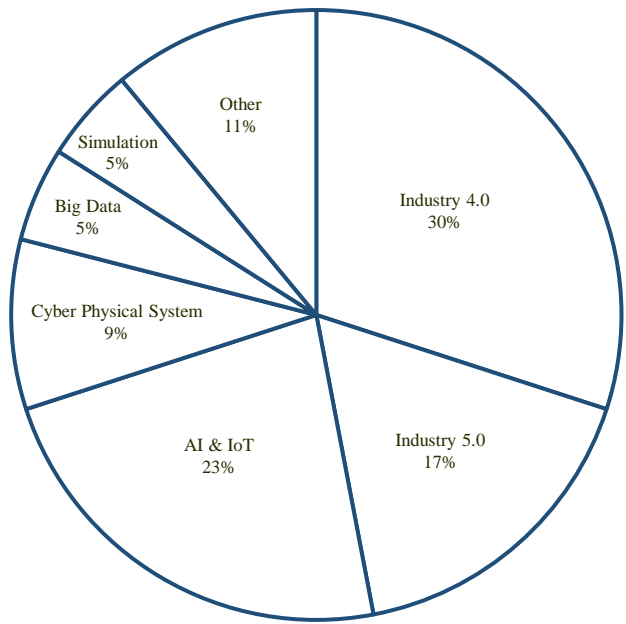

**Figure 10. Use of Technology types of articles.**

## 5. Conclusion

As the topic of SFP has developed and expanded in the new decades of industry, its importance has also been signified. It can be a reason for presenting the framework of the current research. It is critical to shape the framework with the concepts of barriers and enablers according to the I 5.0 background. When the literature is investigated, the presence of smart factories in I 5.0 can be evidence of its importance. This is crucial to forming a successful way to develop the smart factories in I 5.0 because the whole production system is under the meticulous insight of a new generation of technologies. New technologies like AI and IoT and others observe the production system implementation process. If factories tend to have a new way of producing, they must try to apply new technologies in their methods. According to the research, SFP has a nice position in the industry, especially in the future. It is understood that the local predictive automatic response of remote monitoring should be perfectly made in smart factories [66] during the production process. In today's situation, society and people are not too patient to control and monitor things. For this reason, automatic controlling is created. Also, remote monitoring can be productive in the production process.

There are eight enablers and eight barriers in the SFP framework in I 5.0. These factors are extracted according to the articles and the authors' knowledge.

Some enablers include I 5.0 Effects, Production Chain Transformation, Smart Supply, and I 4.0 Technologies. There are some explanations for the enablers here. I 5.0 effects consist of human and machine interaction. The Transformation of the Production Chain consists of modern topics, like facilities and amenities. Sensors implement smart supplies, and digitalization is the key to transforming the supply chain. I 4.0 Technologies can integrate and analyze the instructions in smart factories. Some of the other enablers include data management by processes algorithmizing, moving towards automating new governance models and systems, specialists' skill improvement, and supportive beliefs of environmental society. Data management is an undeniable part of managing roles in modern factories. Also, new governance through self-adaptation and self-optimization is applied. Another enabler is improving the skills of people emerging in new roles in the factories. So, there were some explanations of enablers for readers to be aware of the details.

On the other side, the barriers include Forecasts and decisions monitoring, Challenges in optimizing operations, Lack of environmental resilience, huge investment requirements, Security and privacy, Structural situation instability, Cultural Issues, and Alignment of the amortization process with the business model. Predicting human decisions and preventing acting incorrectly are matters in the SFP. Hardware and software have some requirements that are essential for every smart factory in the industry. Environmental issues and resistance to change are important barriers in the SFP. There are some costs that factories must consider because the cash flow of the companies will be determined. Also, privacy and control are essential for companies because providing a safe environment is crucial for every smart factory. Furthermore, if there is not enough trust among members of factories, then the commitment and mutual understanding will not form properly.

SFP framework as an academic insight has a proper position in the I 5.0. For this reason, some future directions are important to know about it. In the future world, the problem-solving method and the decision-making process will be more precise and accurate according to the studies. Using more modern technologies will make workers and employees more comfortable. So, factories need to be smarter than they are because the world is changing everything; in the same way, workers and humans must be more intelligent than ever in this new era.

As a long term in the Industrial Revolution, the sixth one is coming, and another revolution is waiting for us. Industrial 6.0 (I 6.0) is one of the most important industrial revolutions in the world. Also, Smart factories will tend to use and apply these important generations. The Smart Factory 6.0 can also form a model with concepts. It can probably be an attractive topic for authors. So, it can be stated that the Industrial Revolution and SFP journey did not end here, but it could create a new trend in the whole world. In the future, technology can update more and more, and revolutions and ideas can improve.